\documentclass[10pt]{article}

\usepackage{amsmath,amsfonts,amsthm,commands,graphicx,listings,enumerate,float,multirow}
\usepackage{geometry,algpseudocode,adjustbox,algorithmicx,algorithm,threeparttable,siunitx,color,booktabs}
\usepackage[numbers,sort]{natbib}
\usepackage{subcaption}
\usepackage{lscape}
\usepackage{lipsum}
\usepackage{url}
\newcommand{\I}{\mathcal{I}}
\renewcommand{\L}{\mathcal{L}}

\title{Alternative SDP and SOCP Approximations for Polynomial Optimization}

\author{Xiaolong Kuang
        \thanks{Department of Industrial and Systems Engineering, Lehigh University, Bethlehem, PA, USA.  E-mail: {\tt xik312@lehigh.edu}}
        \and
        Bissan Ghaddar
        \thanks{Department of Management Sciences, University of Waterloo, Canada. E-mail: {\tt bghaddar@uwaterloo.ca}}
        \and
        Joe Naoum-Sawaya
        \thanks{Ivey Business School, University of Western Ontario, Canada. E-mail: {\tt jnaoumsa@uwaterloo.ca}}
        \and
        Luis F. Zuluaga
        \thanks{Department of Industrial and Systems Engineering, Lehigh University, Bethlehem, PA, USA.  E-mail: {\tt luis.zuluaga@lehigh.edu}}
        }

\date{\today}

\begin{document}

% Make title
\maketitle

% Abstract
\begin{abstract}
%A polynomial optimization problem (POP) is an optimization problem whose objective and constraints can be written in terms
%of polynomials on the decision variables.
 In theory,  hierarchies of semidefinite programming (SDP) relaxations  based on sum of squares (SOS) polynomials have been shown to provide arbitrarily close approximations for a general polynomial optimization problem (POP). However, due to the computational challenge of solving SDPs, it becomes difficult to use SDP hierarchies for large-scale problems. % even for low hierarchies' orders.
 To address this, hierarchies of second-order cone programming (SOCP) relaxations resulting from a restriction of the SOS polynomial condition have been recently proposed to approximate POPs. Here, we consider alternative ways to use the SOCP restrictions of the SOS condition. In particular, we show that SOCP hierarchies can be effectively used to strengthen hierarchies of linear programming (LP) relaxations for POPs. Specifically, we show that this solution approach is substantially more effective in finding solutions of certain POPs for which the more common hierarchies of SDP relaxations are known to perform poorly. Furthermore, when the feasible set of the POP is compact, these SOCP hierarchies converge to the POP's optimal value.
\end{abstract}

\smallskip
\noindent \textbf{Keywords\ } Polynomial Optimization, Second-order Cone Relaxation, Semidefinite Relaxation, Approximation Hierarchy.

%*********
% Section
%*********
\section{Introduction}\label{sec.introduction}
Many real-world problems can be modeled as a polynomial optimization problem (POP); that is, an optimization problem in which both the objective and constraints are multivariate polynomials on the decision variables. Thus devising new approaches to globally solve POPs is an active area of research \citep[see, e.g.,][for recent surveys in this area]{BlekPT13, Anjo12}.
In his seminal work, \citet{ref.lasserre} showed that {\em semidefinite programming} (SDP) \cite{Todd01} relaxations based on {\em sum of square} (SOS) polynomials \citep[see, e.g.,][]{BlekPT13} can provide global bounds for POPs. However, the SDP constraints are computationally expensive and thus even using low-orders of the hierarchy to approximate large-scale POPs becomes computationally intractable in practice \cite{lasserre2009moments}. To improve the computational performance of the SDP based hierarchies to approximate the solution of POPs, prior work has focused on exploiting the problem's sparsity \cite{KKW2} and symmetry \cite{DeKlerkS, GatermannP}, improving the relaxation by generating and adding appropriate valid inequalities \cite{ref.bissandigs}, using bounded SOS polynomials \cite{ref.lasserre2015bounded} and more recently by devising more computationally efficient hierarchies such as linear programming (LP) and second-order cone programming (SOCP) hierarchies \cite{ghaddar2011second, ref.zuluaga, ref.amirali, dickinson2013new, dickinson2015extension}.

Here, we consider alternative ways to use SOCP restrictions of the SOS condition introduced by \citet{ref.amirali}. In particular, we show that SOCP hierarchies can be effectively used to strengthen hierarchies of LP relaxations for
 general POPs. Such hierarchies of LP relaxations have received little attention in the POP literature (a few noteworthy
 exceptions are \cite{deKlP02, Lass02e, ZuluVP06,dickinson2013new, dickinson2015extension}). However, in this paper
we show that this solution approach is substantially more effective in finding solutions of certain POPs for which the
 more common hierarchies of SDP relaxations are known to perform poorly
 \citep[see, e.g.,][]{ghaddar2011second}. {Furthermore, when the feasible set of the POP is compact, these SOCP hierarchies converge to the POP's optimal value.} Note that for the well-known SDP based hierarchies introduced by \citet{ref.lasserre}, the {\em quadratic module}~(QM)~\cite{Anjo12} associated with the feasible set of the POP is required to be {\em Archimedean} \cite{BlekPT13}, which implies the compactness of the POP's feasible set.

The remainder of the article is organized as follows. We briefly review several convex approximations of POPs in Section~\ref{sec.litreview}. The proposed approximation strategies and hierarchies are presented in Section~\ref{sec.hierarchy}. Numerical results based on the proposed hierarchies are presented in Section~\ref{sec.results}. Section~\ref{sec.conclusion} provides some concluding remarks.

%*********
% Section
%*********
\section{Preliminaries}\label{sec.litreview}
The following notation is used throughout the article: let $\Rmbb_d[x]:=\Rmbb_d[x_1,\ldots,x_n]$ be the set of polynomials in $n$ variables with real coefficients of degree at most $d$. We define
\bequation\nn
    SOS_{2d}:=\left\{\sum_{i=1}^k p_i(x)^2:p_i(x)\in\Rmbb_{d}[x], k \in\mathbb{Z}_+ \right\},
\eequation
as the cone of SOS polynomials in $\Rmbb_{2d}[x]$ . For any $S\subseteq\Rmbb^n$, let $\Pcal_d(S)$ be the cone of polynomials in $\Rmbb_d[x]$ of degree at most $d$ that are non-negative over the set $S$ \citep[see, e.g.,][]{BlekPT13}. We consider the following general POP,

\bequation\label{opt.ppp}
\tag{PP-P}
    \baligned
         \min_x& & &f(x)  \\
         \st & & &g_i(x)\geq 0,\ i=1,\ldots,m,\\
              & & & x\in\Rmbb^n,
    \ealigned
\eequation
where the degree of the program is $d=\max\{\deg(f),\deg(g_1),\ldots,\deg(g_m)\}$. Given $S=\{x\in\Rmbb^n:g_i(x)\geq0,i=1,\ldots,m\}$, problem~\eqref{opt.ppp} can be equivalently rewritten as the following \emph{conic program} \citep[see, e.g.,][]{ref.convex},

\bequation\label{opt.ppd}
\tag{PP-D}
    \baligned
         \max_{x,\lambda}& & &\lambda  \\
         \st & & &f(x) - \lambda \in \Pcal_d(S),\\
              & & & x\in\Rmbb^n, \lambda\in\Rmbb.
    \ealigned
\eequation
In general, solving~\eqref{opt.ppp} is NP-hard \cite{ref.nphard}. Problem~\eqref{opt.ppd} is a (linear) conic program whose complexity is captured in the cone $\Pcal_d(S)$, which is not {\em tractable} in general. Considering a sequence of {\em tractable} cones $\Kcal^r \subseteq \Kcal^{r+1}\subseteq \cdots \subseteq\Pcal_d(S)$, then the following convex program

\bequation\label{opt.ppm}
%\tag{PP-$\Mcal$}
    \baligned
         \max_{x,\lambda}& & &\lambda  \\
         \st & & &f(x) - \lambda \in \Kcal^r,\\
              & & & x\in\Rmbb^n, \lambda\in\Rmbb
    \ealigned
\eequation
provides a lower bound for~\eqref{opt.ppd}, and hence a lower bound for~\eqref{opt.ppp}. Above, by {\em tractable} we mean that inclusion on the set can be expressed as a {\em linear matrix inequalities} (LMI)~\cite{anjos2012introduction}. As $r$ increases in~\eqref{opt.ppm}, a tighter bound is achieved. The choice of the tractable cone $\Kcal^r$ is a key factor in obtaining good approximation bounds for~\eqref{opt.ppd}, and in turn for~\eqref{opt.ppp}.

%\btheorem[Schm\"{u}dgen~\cite{ref.schm}]
%    Let $f,g_1,...,g_m$ be polynomials in $n$ variables such that $S=\{x\in\mathbb{R}^n:g_i(x)\geq0,i=1,...,m\}$ is compact and $f(x)>0$ for all $x\in S$. Then there exists $s_k(x)\in \Sigma$ such that
%    \bequation\nn
%        f(x)=\sum_k s_k(x)\left(g_1(x)^{\beta_1}\cdots g_m(x)^{\beta_m}\right),
%    \eequation
%    where $\beta_1,\ldots,\beta_m$ are nonnegative integers.
%\etheorem
%Schm\"{u}dgen's theorem implies that we can use tractable sum of square polynomials to approximate the cone $\Pcal_d(S)$.

For this purpose, in his seminal work, \citet{ref.lasserre} proposed a hierarchy of LMI relaxations to approximate $\Pcal_d(S)$, where
\bequation\label{equ.lasscone}
    \Kcal^r =  \left \{p(x) \in \Rmbb_{2r}[x] :p(x)=s_0(x)+\sum_{i=1}^m s_i(x)g_i(x),s_0(x)\in SOS_{2r},s_i(x)\in SOS_{2\lfloor r-\degree(g_i)/2\rfloor}\right \},
\eequation
and $r\geq \lceil d/2 \rceil$ is the level of the hierarchy. In this case, problem~\eqref{opt.ppm} is equivalent to
\bequation\label{opt.lasserre}
\tag{QM-SOS$_r$}
    \baligned
       \max_{x, s_i(x)} & & &\lambda  \\
       \st                  & & &f(x) - \lambda = s_0(x) + \sum_{i=1}^m s_i(x)g_i(x),\\
                             & & &s_0(x)\in SOS_{2r},s_i(x)\in SOS_{2\lfloor r-\degree(g_i)/2\rfloor}, & i=1,\dots, m,\\
                             & & & \lambda\in \Rmbb.
    \ealigned
\eequation
Problem~\eqref{opt.lasserre} can be reformulated as a SDP~\citep[see, e.g.,][]{BlekPT13}. Under some conditions related to the compactness of the set $S$
(more precisely, when the {\em quadratic module} generated by the set of polynomials $\{g_1(x),\ldots,g_m(x)\}$ is {\em Archimedean}), the hierarchy
of problems~\eqref{opt.lasserre}
 converges to the global solution of (\textrm{PP-P}) as $r \to \infty$~\cite{ref.lasserre}. However, as~$r$ increases,
 the size of the positive semidefinite matrices required to reformulate \text{\eqref{opt.lasserre}}  as a SDP increases exponentially.  As a result, this approach is computationally expensive for large-scale problems \citep[see, e.g.,][]{ref.ghaddar2016optimal} or even for small-scale problems that require the solution of high levels of the hierarchy to obtain tight approximations of the POP of interest \citep[see, e.g.,][]{ref.lasserre2015bounded,ref.bissandigs, ref.amirali}.

\citet{ref.amirali} recently proposed a restriction of the SOS condition to address this shortcoming of the SDP-based hierarchies.
The restriction of the SOS condition is done by introducing the use of
%Namely, \citet{ref.amirali} proposed LP- and SOCP-based techniques to approximate sum of square polynomials by introducing
\emph{diagonally dominant sum of square} (DSOS) polynomials and  \emph{scaled diagonally dominant sum of square} (SDSOS) polynomials instead of SOS polynomials in \text{\eqref{opt.lasserre}}.

\begin{definition}[DSOS polynomials~\cite{ref.amirali}]
Let $J$ be an index set, $m_i(x) \in \Rmbb_d[x]$ be a monomial for all $i\in J$, and
$\alpha_i, \beta_{ij} \in \Rmbb_+$ for all $i,j \in J$. Then
\begin{equation}
\label{eq:DSOS}
p(x) = \sum_i\alpha_i m_i(x)^2 + \sum_{i,j}\beta_{ij}(m_i(x) \pm m_j(x))^2,
\end{equation}
is a DSOS polynomial in $\Rmbb_{2d}[x]$.
Equivalently, DSOS polynomials can be defined as those that can be constructed from a {\em diagonally dominant
matrix} (DD). Namely, let $z(x)$ be a vector with the monomials $m_i(x)$ for all $i \in J$, and $Q \in \Rmbb^{|J| \times |J|}$
be a (symmetric) diagonally dominant matrix. Then $p(x)=z^T(x)Qz(x)$ is a DSOS.
\end{definition}

Let $DSOS_{2d}$ be the set of all DSOS polynomials in $\Rmbb_{2d}[x]$. Then it is clear from \eqref{eq:DSOS} that $DSOS_{2d} \subseteq SOS_{2d}$.
Thus, using DSOS polynomials instead of SOS polynomials in \text{\eqref{opt.lasserre}} provides a hierarchy of
lower bounds for the SOS hierarchy. Moreover the resulting DSOS hierarchy is computationally easier to solve. Namely, recall
 that a symmetric matrix $A\in \Rmbb^{n \times n}$ is DD if $A_{ii}\geq\sum_{j\neq i}|A_{ij}|,\forall i=1,\dots,n$.
 Thus the associated DSOS hierarchy
 \bequation\label{opt.lp}
\tag{QM-DSOS$_r$}
    \baligned
         \max_{\lambda, d_i(x)}& & &\lambda  \\
         \st                  & & &f(x) - \lambda = d_0(x) + \sum_{i=1}^m d_i(x)g_i(x),\\
                     & & &d_0(x)\in DSOS_{2r},d_i(x)\in DSOS_{2\lfloor r-\degree(g_i)/2\rfloor},\\
                             & & & \lambda\in \Rmbb,
    \ealigned
\eequation
can be reformulated as a LP. As proposed by \citet{ref.amirali}, the DSOS hierarchy~\eqref{opt.lp} can be strengthened
by considering \emph{scaled diagonally dominant sum of square} (SDSOS) polynomials.

\begin{definition}[SDSOS polynomials~\cite{ref.amirali}]
Let $J$ be an index set, $m_i(x) \in \Rmbb_d[x]$ be a monomial for all $i\in J$, and
$\alpha_i, \beta_i, \beta_{j} \in \Rmbb_+$ for all $i,j \in J$. Then
\begin{equation}
\label{eq:SDSOS}
p(x) = \sum_i\alpha_i m_i(x)^2 + \sum_{i,j}(\beta_{i}m_i(x) \pm \beta_{j}m_j(x))^2,
\end{equation}
is a SDSOS polynomial in $\Rmbb_{2d}[x]$.
Equivalently, SDSOS polynomials can be defined as those that can be constructed from a {\em scaled diagonally dominant
matrix} (SDD). Namely, let $z(x)$ be a vector with the monomials $m_i(x)$ for all $i \in J$, and $Q \in \Rmbb^{|J| \times |J|}$
be a (symmetric) scaled diagonally dominant matrix. Then $p(x)=z^T(x)Qz(x)$ is a SDSOS.
\end{definition}

Let $SDSOS_{2d}$ be the set of all SDSOS polynomial in $\Rmbb_{2d}[x]$. Then it is clear from
\eqref{eq:SDSOS} that $DSOS_{2d} \subseteq SDSOS_{2d} \subseteq SOS_{2d}$.
Thus, using SDSOS polynomials instead of SOS polynomials in \text{\eqref{opt.lasserre}} provides a hierarchy of
lower bounds for the SOS hierarchy that is tighter than the~\eqref{opt.lp} hierarchy. Moreover the resulting SDSOS hierarchy is computationally easier to solve than the associated SDP-based hierarchy. Namely, notice that
 a symmetric matrix $A\in \Rmbb^{n \times n}$ is SDD if
 \begin{equation}
 \label{eq:SDSOSmat}
 A = \sum_{i, j \in \{1,\dots,n\}} A^{ij},  \text{ for some } A^{ij} \succeq 0, \text{ with } A^{ij}_{kl} = 0  \text{ for
 any } k,l \in \{1,\dots,n\} \setminus \{i,j\}.
 \end{equation}
 Above, we use the common notation $A \succeq 0$ to indicate that the matrix is positive semidefinite. Notice that because the $A^{ij}$ matrices in~\eqref{eq:SDSOSmat} have only nonzero elements
 at positions $k,l \in \{i,j\}$, then it follows that

 \bequation\label{equ.sdsostosocp}
        A^{ij} \succeq 0  \iff   A^{ij}_{ii}+A^{ij}_{jj}\geq \left\|\begin{pmatrix}2A^{ij}_{ij} \\ \\
                                                                                    A^{ij}_{ii}-A^{ij}_{jj} \end{pmatrix}\right\|_2
                                                                                     \iff \begin{pmatrix} A^{ij}_{ii}+A^{ij}_{jj} \\
                                                                                    2A^{ij}_{ij}\\
                                                                                    A^{ij}_{ii}-A^{ij}_{jj} \end{pmatrix} \in \L^3,
\eequation
where $\L^n$ denotes the {\em second-order cone} or {\em Lorentz cone} of dimension $n$~\citep[see, e.g.,][]{ref.convex}. Thus the associated SDSOS hierarchy
 \bequation
 \label{opt.socp}
\tag{QM-SDSOS$_r$}
    \baligned
        \max_{\lambda, d_i(x)}& & &\lambda  \\
        \st                  & & &f(x) - \lambda = d_0(x) + \sum_{i=1}^m d_i(x)g_i(x),\\
                    & & &d_0(x)\in SDSOS_{2r},d_i(x)\in SDSOS_{2\lfloor r-\degree(g_i)/2\rfloor},\\
                    & & & \lambda\in\Rmbb,
    \ealigned
\eequation
can be reformulated as a second-order cone program. \citet{ref.amirali} have shown that the approximation hierarchies \eqref{opt.lp} and \eqref{opt.socp} can be successfully used to approximate POPs arising in control, combinatorics, and general non-linear non-convex optimization~\cite{ref.amirali}. Hierarchies~\eqref{opt.lp} and~\eqref{opt.socp} are computationally easier to solve than~\eqref{opt.lasserre}, however, their bounds might not be as good as the one obtained with the~\eqref{opt.lasserre} hierarchy of the same order~\citep[see, e.g.,][]{ref.kuang2016alternative}.

%*********
% Section
%*********

%*********
% Section
%*********
\section{Alternative LP, SOCP and SDP Hierarchies for POP}\label{sec.hierarchy}
Lasserre's hierarchy~\cite{ref.lasserre} has been shown to provide very tight bounds
for multiple classes of POPs.
However, this approach becomes computationally intractable for large-scale problems
or even for small-scale problems that require the solution of high levels of the hierarchy to obtain good approximations for the solution of the problem of interest. Loosely speaking, this intractability stems from the fact that the size of the SDP reformulation of the SOS conditions in  \text{\eqref{opt.lasserre}}
grows exponentially with the dimension of the decision variables of the problem~$n$, as well as
the level of the hierarchy~$r$.

A key building block behind the  convergence properties of the hierarchy defined by \text{\eqref{opt.lasserre}} is a representation
theorem for polynomials in $\Pcal_d(S)$ by \citet{Puti93} that makes use of SOS polynomials~\citep[see, e.g.,][]{BlekPT13, ref.lasserre}. Other
convergent SDP hierarchies can be constructed similarly using the representation theorem by \citet{Schm91},  when the
set $S$ is compact. Besides these SOS representation theorems, there are however well-known representations theorems for
non-negative polynomials that use polynomials with non-negative coefficients (instead of SOS polynomials) in the representation. Examples of these are the representation theorem of \citet{Poly88}, when the set $S$ is a polytope, and P{\'o}lya's Theorem~\cite{Poly88}, when the set $S = \Rmbb^n_+$.

\begin{theorem}[P\'olya \citep{Poly88}]
\label{thm:polya}
Let $p(x) \in \Rmbb^n[x]$ be a multivariate polynomial. Then
$p(x) > 0$ \text{ for all } $x \ge 0 \Rightarrow \left (1 + \sum_{i=1}^n x_i \right )^r p(x) = \sum_{\alpha \in \Nmbb^n} c_{\alpha} x^{\alpha}$	
for some $r \ge 0$, $c_{\alpha} \ge 0$ for all $\alpha \in \Nmbb^n$.
\end{theorem}

In stating Theorem~~\ref{thm:polya}, we make use of the common notation
$x^\alpha := x_1^{\alpha_1} \cdots x_n^{\alpha_n}$ for any $x \in \Rmbb^n$, and
$\alpha \in \Nmbb^n$. Note that in Theorem~\ref{thm:polya},
the non-negativity of the polynomial is certified using
polynomials with non-negative coefficients. As a result, this type of representation theorems can be used to construct hierarchies of LP problems that converge to the optimal solution of~\eqref{opt.ppp} (when the required conditions on the set $S$ are satisfied). Such approach has been used in~\cite{deKlP02, Lass02e, ZuluVP06}. It is worthy to mention that P\'olya's approach is also used in~\cite{dickinson2015extension,dickinson2013new}, to address the solution of POPs.

Here, we take advantage of this type of computationally easier LP hierarchy approach to address the solution of certain classes of POPs for which
the more common SDP hierarchy is known to perform poorly~\citep[see, e.g.,][]{ghaddar2011second}.
In particular, we use a representation theorem for non-negative polynomials in a semi-algebraic set
recently  introduced in~\cite{ref.zuluaga} to construct
a converging hierarchy of LPs for POPs. Formally, consider the following optimization problem:
\bequation\label{opt.zuluaga}
\tag{Po-LP$_r$}
    \baligned
           z_{r,LP}:=\max_{\lambda, c_{\alpha, \beta}}& & &\lambda \\
               \st & & &\left (1 + \sum_{i=1}^n x_i+\sum_{j=1}^m g_j(x) \right )^r (f(x) - \lambda) = \sum_{(\alpha,\beta)\in \I} c_{\alpha,\beta}x^\alpha g(x)^{\beta}\\
                   & & &c_{\alpha,\beta} \in\Rmbb_+ \text{ for all } (\alpha,\beta)\in \I, \\
                   & & & \lambda\in\Rmbb,
              %& &     & & & t=\left\lfloor \frac{r-\degree(f)}{\max(\degree(g))} \right\rfloor,
    \ealigned
\eequation
where
\[
\I : = \{(\alpha,\beta)\in\Nmbb^{n+m}: \|(\alpha,\beta)\|_1 \leq r\max\{\deg(g_i): i =1,\ldots,m\} + \deg(f)\}.
\]
By matching the coefficients of each monomial in the left-hand side and the right-hand side of equation~\eqref{opt.zuluaga}, the resulting problem is a LP with decision variables $\lambda \in \Rmbb, c_{\alpha,\beta} \in \Rmbb_+$, for all $(\alpha,\beta)\in\I$.
Similar to the  \text{\eqref{opt.lasserre}} hierarchy, but under milder conditions, the resulting bound~$z_{r,LP}$ of~\eqref{opt.zuluaga}, obtained at each level of the hierarchy converges as~$r$ increases. \btheorem[\citet{ref.zuluaga}]\label{thm.zuluaga}
    Let $S=\{x\in \Rmbb^n_+:g_i(x)\geq 0,i=1,\ldots,m\}$ be a compact set, then as $r\rightarrow\infty$, $ z_{r,LP}$ converges to the global optimum of~\eqref{opt.ppp}.
\etheorem

Thus, Theorem~\ref{thm.zuluaga} provides the convergence guarantee of
the~\eqref{opt.zuluaga} hierarchy to the optimal value of (PP-P) with a compact feasible set in~$\Rmbb^+_n$. This allows us to use LP techniques to globally solve non-convex problems. However, this type of LP approximations
for POPs are known to provide very weak approximation bounds for the objective value of the
POP of interest \citep[see, e.g.,][]{deKlP02, Lass02e}. To address this, we next propose
the use of DSOS, SDSOS and SOS polynomials with fixed degree (degree 2) instead of the non-negative constant
$c_{\alpha,\beta}$ in the definition of the hierarchy \eqref{opt.zuluaga}.

For a general POP~\eqref{opt.ppp} with feasible set $ S \subseteq \Rmbb^n_+$, consider the
following hierarchies of optimization problems:
\bequation\label{opt.sos2}
\tag{Po-SOS$_r$}
    \baligned
         \max_{\lambda, p_{\alpha, \beta}}& & & \lambda \\
              \st &  & & \left (1 + \sum_{i=1}^n x_i+\sum_{j=1}^m g_j(x) \right )^r (f(x) - \lambda) = \sum_{(\alpha,\beta)\in \I'} p_{\alpha,\beta}(x) x^\alpha g(x)^{\beta},\\
             &     & & p_{\alpha,\beta}(x)\in SOS_2, \text{ for all } (\alpha,\beta)\in \I',\\
             & &  & \lambda\in\Rmbb,
    \ealigned
\eequation

\bequation\label{opt.sdsos2}
\tag{Po-SDSOS$_r$}
    \baligned
         \max_{\lambda, p_{\alpha, \beta}}&  & & \lambda \\
              \st & & & \left (1 + \sum_{i=1}^n x_i+\sum_{j=1}^m g_j(x) \right )^r (f(x) - \lambda) = \sum_{(\alpha,\beta)\in \I'} p_{\alpha,\beta}(x) x^\alpha g(x)^{\beta},\\
                  & & & p_{\alpha,\beta}(x)\in SDSOS_2, \text{ for all } (\alpha,\beta)\in \I',\\
             & &  & \lambda\in\Rmbb,
    \ealigned
\eequation

\bequation\label{opt.dsos2}
\tag{Po-DSOS$_r$}
    \baligned
         \max_{\lambda, p_{\alpha, \beta}}&  & & \lambda \\
              \st &  & & \left (1 + \sum_{i=1}^n x_i+\sum_{j=1}^m g_j(x) \right )^r (f(x) - \lambda) = \sum_{(\alpha,\beta)\in \I'} p_{\alpha,\beta}(x) x^\alpha g(x)^{\beta},\\
             &     & & p_{\alpha,\beta}(x)\in DSOS_2, \text{ for all } (\alpha,\beta)\in \I',\\
             & &  & \lambda\in\Rmbb,
    \ealigned
\eequation

where $r\geq0$ and
\[
\I' : = \{(\alpha,\beta)\in\Nmbb^{n+m}: \|(\alpha,\beta)\|_1 \leq r\max\{\deg(g_i): i=1,\ldots,m\} + \deg(f)-2\}.
\]
Similar to Lasserre's hierarchy~\eqref{opt.lasserre}, problem~\eqref{opt.sos2} can be reformulated as a SDP. In turn, similar to the hierarchies~\eqref{opt.lp} and~\eqref{opt.socp} (cf., Section~\ref{sec.litreview}),
the optimization problems~\eqref{opt.dsos2} and~\eqref{opt.sdsos2} can be reformulated as a LP and as a SOCP respectively.
Note that in the hierarchies discussed in Section~\ref{sec.litreview}, as the level of the
hierarchy $r$ increases, the complexity of checking that a fixed number, $m+1$, of polynomials are SOS, SDSOS, or DSOS increases. Instead in the hierarchy defined
in~\eqref{opt.sos2},~\eqref{opt.sdsos2} and~\eqref{opt.dsos2}, the complexity of checking that the involved polynomials are SOS, SDSOS, or DSOS does not change as the degree of these polynomials is fixed to 2. Instead, it is the number of these
polynomials that increases as the level of the hierarchy increases (a similar approach has been used in~\cite{ref.lasserre2015bounded}). This turns out to be key to obtain the results presented later in next section on the performance of the hierarchies~\eqref{opt.sos2},~\eqref{opt.sdsos2} and~\eqref{opt.dsos2}.

Clearly, the hierarchies~\eqref{opt.sos2},~\eqref{opt.sdsos2}, and~\eqref{opt.dsos2} provide tighter bounds on the associated POP than
the LP based hierarchy~\eqref{opt.zuluaga}. As a result, under the same conditions of Theorem~\ref{thm.zuluaga},
these hierarchies will converge as $r \to \infty$ to the global optimal solution of~\eqref{opt.ppp}.
Below, we state this formally.

\bproposition\label{pro.bounds}
Consider problem~\eqref{opt.ppp} with a compact feasible region and assume that $S \subseteq \Rmbb^n_+$
 whose global optimal objective function is $z^\ast$, and let $ z_{r,DSOS}$, $z_{r,SDSOS}$, $z_{r,SOS}$ be the optimal value of hierarchies~\eqref{opt.dsos2},~\eqref{opt.sdsos2} and~\eqref{opt.sos2} respectively, then it follows that for any $r=1,2,\dots$:
    \bequation\nn
        z_{r,LP}\leq z_{r,DSOS} \leq z_{r,SDSOS} \leq z_{r,SOS}\leq z^\ast.
    \eequation
    Moreover,
  \[
     \lim_{r \to \infty}  z_{r,DSOS} = \lim_{r \to \infty}  z_{r,SDSOS} = \lim_{r \to \infty}  z_{r,SOS} =  z^\ast.
    \]
\eproposition
\bproof
The inequalities $z_{r,DSOS} \leq z_{r,SDSOS} \leq z_{r,SOS}$ follow from $DSOS_{2d}\subseteq SDSOS_{2d}\subseteq SOS_{2d}$. It is easy to see $z_{r,LP}\leq z_{r,DSOS}$ since all the nonnegative constants belong to $DSOS_{0}$. By Theorem~\ref{thm.zuluaga}, $\lim_{r\rightarrow\infty} z_{r,LP}=z^\ast$ when the feasible region of~\eqref{opt.ppp} is compact. Thus $\lim_{r \to \infty}  z_{r,DSOS} = \lim_{r \to \infty}  z_{r,SDSOS} = \lim_{r \to \infty}  z_{r,SOS} =  z^\ast$ follows.
\eproof
\section{Numerical Results}\label{sec.results}
To illustrate the performance of the hierarchies discussed in Section~\ref{sec.hierarchy}, we test the Lasserre-type hierarchies and the proposed hierarchies in this article on some relevant POP instances. We use
APPS~\cite{ref.bissanapps} together with {\tt Matlab} to implement all the hierarchies. Numerical experiments are conducted on an AMD Opteron 2.0 GHz(x16) Linux machine with 32 GB memory. We use MOSEK~\cite{ref.mosek} as
the LP and SOCP solver. Also, we use SeDuMi~\cite{ref.sedumi} as the SDP solver, to exploit its well-known precision for solving SDPs.

Due to the different approach used in the Lasserre-type hierarchies and the hierarchies proposed in Section~\ref{sec.hierarchy}, with the same $r$, the degree of the polynomials involved in the problem might not be equal. Thus, to make it easier to compare the results obtained from each hierarchy, instead of reporting the hierarchy level $r$, we report the maximum degree $\hat{d}$ of the polynomials involved in the formulation as~$r$ increases in each of them. %To address this imbalnace, we make a little change in running the experiments. Recall $d$ as the maximum degree of objective and constraints, we start with all the hierarchies with degree $\hat{d}:=2\lceil d/2\rceil$ and assign $\hat{d}:=\hat{d}+2$ for next level. With the increase of $\hat{d}$, increase $r$ in both types of hierarchies accordingly.

In the tables of numerical results that follow, the symbol (*) indicates that the reported value is the optimal objective value of the problem. We use ``T'' as the solution time in seconds for each hierarchy and ``Infeas.'' to indicate that the optimization problem is infeasible. The symbol ($\diamond$) indicates that the solver runs out of memory. Lastly the symbol ($\circ$) indicates that generating the program that matches coefficient in the hierarchy in {\tt Matlab} runs out of memory.
\subsection{Illustrative Examples}\label{sec.numsub1}
We begin by testing a set of POPs from~\cite{ref.bissandigs}, which are highly non-convex and require a high level of Lasserre's hierarchy to converge to their global optimum.
\begin{example}
Consider the following quadratic POP with 5 variables:
\bequation\nn %\label{exa.5vars}
    \begin{split}
        \min_{x \in {\Rmbb}^5}\ \ & 2x_1-x_2+x_3-2x_4-x_5\\
        \st   \ \ & (x_1-2)^2-x_2^2-(x_3-1)^2-(x_5-1)^2\geq0,\\
                  & x_1x_3-x_4x_5+x_1^2\geq1,\\
                  & x_3-x_2^2-x_4^2\geq1,\\
                  & x_1x_5-x_2x_3\geq2,\\
                  & x_1+x_2+x_3+x_4+x_5\leq14,\\
                  & x_i\geq0,i=1,\ldots,5.
    \end{split}
\eequation
\end{example}
As shown in Table~\ref{tab.degree5vars}, the~\eqref{opt.lasserre} hierarchy converges to the global optimum  when $\hat{d}=8$ with a computational time of 49.82 seconds, while the hierarchy~\eqref{opt.sos2} converges to global optimum when $\hat{d}=6$ with only 8.21 seconds of computational time. Hierarchies~\eqref{opt.socp} and~\eqref{opt.lp} fail to converge to the global optimum up to $\hat{d}=8$. However, the hierarchy~\eqref{opt.sdsos2} also converges to the global optimum when $\hat{d}=8$ with 13.28 seconds of computational time. The hierarchy~\eqref{opt.dsos2} provides a weaker bound than hierarchy~\eqref{opt.sdsos2} and does not converge to the problem's global optimum when $\hat{d}=8$.

\begin{figure}[!tb]
  \begin{subfigure}[b]{0.5\linewidth}
    \centering
    \includegraphics[width=\linewidth]{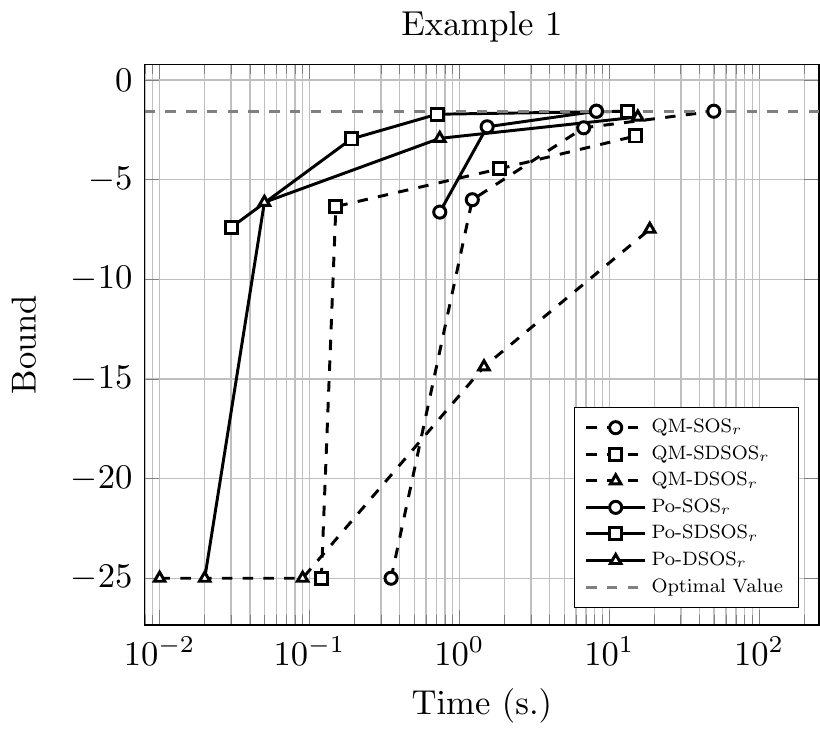}
   % \caption{Initial condition}
    \label{fig.compare:a}
  \end{subfigure}%%
  \begin{subfigure}[b]{0.5\linewidth}
    \centering
    \includegraphics[width=\linewidth]{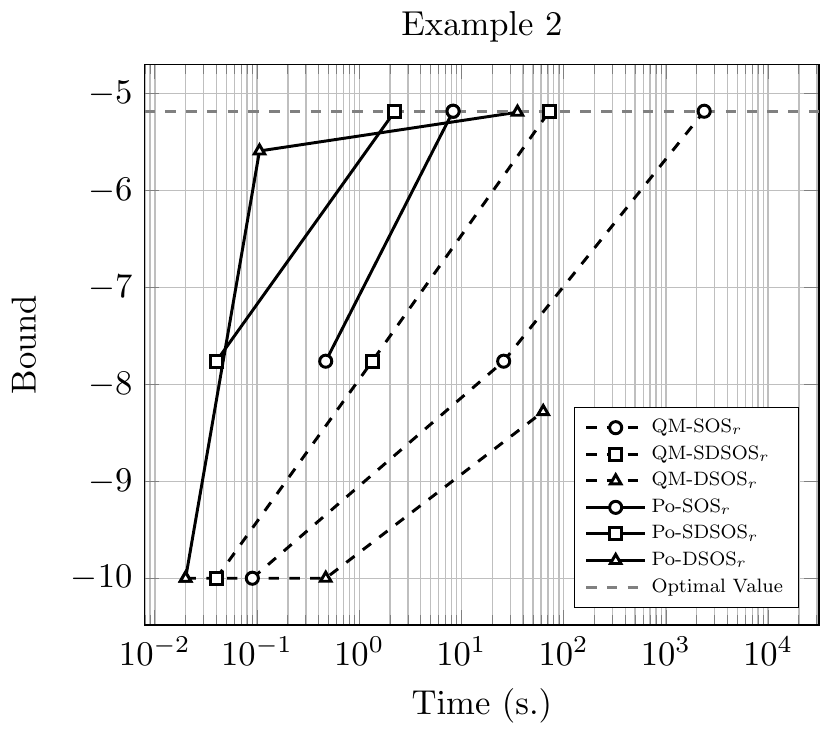}
  %  \caption{Rupture}
    \label{fig.compare:b}
  \end{subfigure}
  \caption{Bound and Time Comparison of Different Hierarchies for Example 1 (left) and Example~2 (right).}
  \label{fig.compare}
\end{figure}

Although the degree $\hat{d}$ provides an approximate measure of the size (variables and constraints) involved in the formulations of the hierarchies' problems, a better comparison of the hierarchies can be done by illustrating the trade-off between the solution time and the quality of the bound obtained from each hierarchy. In
Figure~\ref{fig.compare} (left), the different line plots show the bound and solution time associated with increasing orders of each of the hierarchies. Clearly, within one second, the~\eqref{opt.sdsos2} hierarchy gives the best bound; within ten seconds, the~\eqref{opt.sos2} hierarchy gives the optimal value while there is still a gap between the problem's optimal value
(illustrated by the dashed horizontal line) and the bounds obtained by other hierarchies. Clearly, the hierarchies proposed in Section~\ref{sec.hierarchy} have better performance over the Lasserre-type hierarchies for this problem.

%\begin{table}[H]
%    \centering
%
%    \begin{tabular}{lrrrrrr}
%        \toprule
%        &\multicolumn{2}{c}{SDP hierarchy}&\multicolumn{2}{c}{SOCP hierarchy}&\multicolumn{2}{c}{LP hierarchy}\\
%        r & Bound & T & Bound & T & Bound & T \\\midrule
%        2 & 25.00 & 0.34 & 25.00 & 0.09 &25.00 & 0.03   \\
%        4 & 6.01 & 1.13 & 6.66 & 0.19 & 25.00 & 0.11  \\
%        6 & 2.40 & 5.43 & 4.55 & 2.13 & 14.39 & 2.09  \\
%        8 & $^\ast$1.57 & 37.32 & 3.07 & 17.81 & 7.49 & 26.34   \\
%        10 &            &       &  2.06 & 129.38 & 4.06& 394.90 \\
%        \bottomrule
%    \end{tabular}
%    \caption{Comparison of Lasserre's Hierarchy and LP-based, SOCP-based Hierarchies for Example 1.}\label{tab.lpsocp5vars}
%\end{table}

\begin{table}[H]
    \centering
   \begin{adjustbox}{width=1\textwidth}
    \begin{tabular}{crrrrrrrrrrrr}
        \toprule
         &\multicolumn{2}{c}{\eqref{opt.lasserre}}&\multicolumn{2}{c}{\eqref{opt.socp}}&\multicolumn{2}{c}{\eqref{opt.lp}}&\multicolumn{2}{c}{\eqref{opt.sos2}}&\multicolumn{2}{c}{\eqref{opt.sdsos2}}&\multicolumn{2}{c}{\eqref{opt.dsos2}}\\ \midrule
         $\hat{d}$ & Bound & T & Bound & T  & Bound & T & Bound & T & Bound & T & Bound & T\\\midrule
          2 & -25.00 & 0.35 & -25.00 & 0.12 & -25.00 & 0.01 & -6.63 & 0.74 & -7.40 & 0.03 & -25.00 & 0.02 \\
         4 & -6.01 & 1.22 & -6.35 & 0.15 & -25.00 & 0.09 & -2.35 & 1.53 & -2.96 & 0.19 & -6.14 & 0.05  \\
         6 & -2.40 & 6.75 & -4.46 & 1.85 & -14.39 & 1.46 & $^\ast$-1.57 & 8.21 & -1.72 & 0.71 & -2.93 & 0.74  \\
         8 & $^\ast$-1.57 & 49.82 & -2.81 & 15.00 & -7.49 & 18.62 &  &  & $^\ast$-1.57 & 13.28 & -1.86 & 15.49   \\
        \bottomrule
    \end{tabular}
  \end{adjustbox}
\begin{tablenotes}
        \footnotesize
        \item[1] $\ast$: Optimal value is obtained.
    \end{tablenotes}
\caption{Bound and Time Comparison of Different Hierarchies for Example 1.}\label{tab.degree5vars}
\end{table}
In Table~\ref{tab.degree5vars}, note that for the same level of hierarchies~\eqref{opt.dsos2} and~\eqref{opt.sdsos2}, the linear representation of $DSOS_{2}$ introduces more decision variables than the SOCP representation of $SDSOS_{2}$. This explains why the running time of the LP-based hierarchy can be larger than the running time of the SOCP-based hierarchy.

\begin{example}
Consider the following quadratic POP with 10 variables:
\bequation\nn %\label{exa.10vars}
    \begin{split}
        \min_{x \in {\Rmbb^{10}}}\ \ & -x_1-x_2+x_3-2x_4-x_5+x_6+x_7-x_8+x_9-2x_{10}\\
        \st   \ \ & (x_3-2)^2-(x_5-1)^2-2x_6+x_8^2-(x_9-2)^2\geq-4,\\
                  & -x^2_2+x_3x_{10}-x_4^2+x_6x_7\geq1,\\
                  & x_1x_8-x_2x_3+x_4x_7-x_5x_{10}\geq2,\\
                  & \sum_{i=1}^{10} x_i\leq5,\\
                  & x_i\geq0,i=1,\ldots,10.
    \end{split}
\eequation
\end{example}
As shown in Table~\ref{tab.degree10vars}, Lasserre's hierarchy~\eqref{opt.lasserre} and hierarchy~\eqref{opt.socp} converge to the global optimum at the third level when $\hat{d}=6$ with a computational time of 2369.50 seconds and 72.43 seconds respectively. In contrast, hierarchy~\eqref{opt.sos2} converges to the global optimum
when $\hat{d}=4$ with 8.27 seconds of computational time. The hierarchy~\eqref{opt.sdsos2} also converges to the global optimum when $\hat{d}=4$ with 2.23 seconds of computational time.
Similar to Example 1, hierarchies~\eqref{opt.lp} and~\eqref{opt.dsos2} provide the weakest bound and the problem's global optimum is not reached by $\hat{d}=6$, but the hierarchy~\eqref{opt.dsos2} provides tighter bounds with less computational time than the hierarchy~\eqref{opt.lp} at each level.

As discussed previously, a better comparison among the different hierarchies can be obtained by
illustrating the trade-off between the solution time and the quality of the bound obtained from each hierarchy.
In Figure~\ref{fig.compare} (right),
the different line plots show the bound  and solution time associated with increasing orders of each of the hierarchies. Notice that
within one second, the~\eqref{opt.dsos2} gives the best bound. Also, within ten seconds, only the~\eqref{opt.sos2} and~\eqref{opt.sdsos2} hierarchies obtain the problem's optimal value (illustrated by the dashed horizontal line), and the~\eqref{opt.sdsos2} hierarchy takes less computational time than the~\eqref{opt.sos2} hierarchy.
\\

\begin{table}[!tb]
\centering
   \begin{adjustbox}{width=1\textwidth}
    \begin{tabular}{crrrrrrrrrrrr}
        \toprule
         &\multicolumn{2}{c}{\eqref{opt.lasserre}}&\multicolumn{2}{c}{\eqref{opt.socp}}&\multicolumn{2}{c}{\eqref{opt.lp}}&\multicolumn{2}{c}{\eqref{opt.sos2}}&\multicolumn{2}{c}{\eqref{opt.sdsos2}}&\multicolumn{2}{c}{\eqref{opt.dsos2}}\\ \midrule
         $\hat{d}$ & Bound & T & Bound & T  & Bound & T & Bound & T & Bound & T & Bound & T\\\midrule
          2 & -10.00 & 0.09 & -10.00 & 0.04 & -10.00 & 0.02 & -7.76 & 0.17 & -7.76 & 0.04 & -10.00 & 0.02 \\
         4 & -7.76 & 25.89 & -7.76 & 1.34 & -10.00 & 0.47  & $^\ast$-5.18 & 8.27 & $^\ast$-5.18 & 2.23 & -5.59 & 0.16  \\
         6 & $^\ast$-5.18 & 2369.50 & $^\ast$-5.18 & 72.43 & -8.28 & 63.09 &  &  &  &  & -5.19 & 35.29 \\
        \bottomrule
    \end{tabular}
  \end{adjustbox}
      \begin{tablenotes}
        \footnotesize
        \item[1] $\ast$: Optimal value is obtained.
    \end{tablenotes}
\caption{Bound and Time Comparison of Different Hierarchies for Example 2.}\label{tab.degree10vars}
\end{table}

\begin{example}
Consider the following quadratic POP with 15 variables:
\bequation\nn %\label{exa.15vars}
    \begin{split}
        \min_{x \in \Rmbb^{15}}\ \ & x_1-x_2+x_3-x_4-x_5+x_6+x_7-x_8+x_9-x_{10}+x_{11}-x_{12}+x_{13}-x_{14}+x_{15}\\
        \st   \ \ & (x_1-2)^2-x_2^2+(x_3-2)^2-(x_4-1)^2-(x_5-1)^2+(x_6-2)^2-(x_7-1)^2-x_8^2\\
                  & -(x_9-2)^2-(x_{10}-1)^2+x_{11}^2-x_{12}^2+(x_{13}-2)^2+x_{14}^2-(x_{15}-1)^2\geq0,\\
                  & -x_1x_7-x_4x_5-x_3^2+x_6x_9+x_{10}x_{12}\geq1,\\
                  & x_2x_3-x_8x_{11}-x_{14}^2+x_5x_{15}\geq2,\\
                  & \sum_{i=1}^{15} x_i\leq10,\\
                  & x_i\geq0,i=1,\ldots,15.
    \end{split}
\eequation
\end{example}

The results for Example 3 are shown in Table \ref{tab.degree15vars}. Lasserre's hierarchy~\eqref{opt.lasserre} and the hierarchy~\eqref{opt.socp} fail to provide the problem's global optmial value when $\hat{d}=4$. {\tt Matlab} runs out of memory when generating the LMI for Lasserre-type hierarchies when $\hat{d}=6$. In contrast, hierarchies~\eqref{opt.sos2} and~\eqref{opt.sdsos2} converge to the global optimum when $\hat{d}=4$ with 640.60 and 59.85 seconds of computational time respectively. Similar to Example 1 and 2, hierarchies~\eqref{opt.lp} and~\eqref{opt.dsos2} provide the weakest bound and the problem's global optimum is not reached when $\hat{d}=6$. However, the~\eqref{opt.dsos2} hierarchy provides tighter bounds with less computational time than~\eqref{opt.lp} when $\hat{d}=2$ and $\hat{d}=4$.  \\

\begin{table}[!tb]
    \centering
   \begin{adjustbox}{width=1\textwidth}
    \begin{tabular}{crrrrrrrrrrrr}
        \toprule
         &\multicolumn{2}{c}{\eqref{opt.lasserre}}&\multicolumn{2}{c}{\eqref{opt.socp}}&\multicolumn{2}{c}{\eqref{opt.lp}}&\multicolumn{2}{c}{\eqref{opt.sos2}}&\multicolumn{2}{c}{\eqref{opt.sdsos2}}&\multicolumn{2}{c}{\eqref{opt.dsos2}}\\ \midrule
         $\hat{d}$ & Bound & T & Bound & T  & Bound & T & Bound & T & Bound & T & Bound & T\\\midrule
          2 & -10.00 & 0.09 & -10.00 & 0.02 & -10.00 & 0.02 & -8.07 & 0.27  & -8.74 & 0.06 & -10.00 & 0.02 \\
         4 & -8.06 & 2754.30 & -8.29 & 10.17 & -10.00 & 2.27 & $^\ast$-7.43 & 640.60 & $^\ast$-7.43 & 59.85 & -8.22 & 0.51  \\
         6 &$\circ$ &$\circ$ &$\circ$  & $\circ$ & $\circ$ & $\circ$ &  &  &  & & -7.64 & 2340.00  \\
        \bottomrule
    \end{tabular}
  \end{adjustbox}
    \begin{tablenotes}
        \footnotesize
         \item[1] $\ast$: Optimal value is obtained.
        \item[2] $\circ$: {\tt Matlab} runs out of memory while formulating LMI.
    \end{tablenotes}
\caption{Bound and Time Comparison of Different Hierarchies for Example 3.}\label{tab.degree15vars}

\end{table}

\subsection{Numerical Results on Global Optimization Library}\label{sec.numsub2}
Next, we compare Lasserre-type hierarchies with the proposed hierarchies on some problems from the GLOBAL Library available at \url{http://www.gamsworld.org/global/globallib.htm}. These problems have been used as benchmark in~\cite{waki2006sums, waki2008algorithm, kleniati2010partitioning}.

\begin{figure}[!tb]
  \begin{subfigure}[b]{0.5\linewidth}
    \centering
    \includegraphics[width=\linewidth]{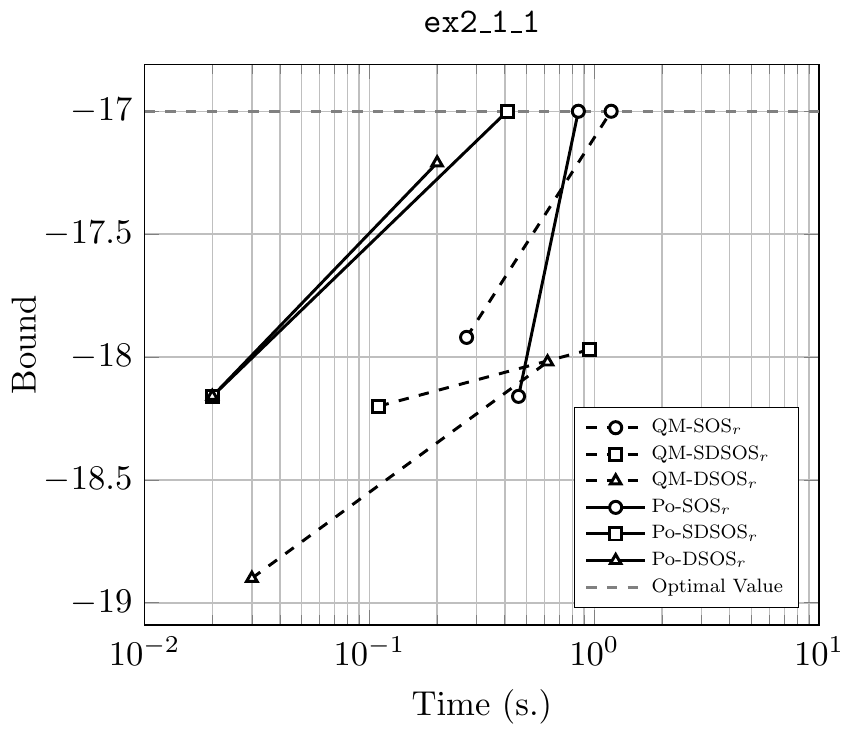}
   % \caption{Initial condition}
    \label{fig7:a}
  \end{subfigure}%%
  \begin{subfigure}[b]{0.5\linewidth}
    \centering
    \includegraphics[width=\linewidth]{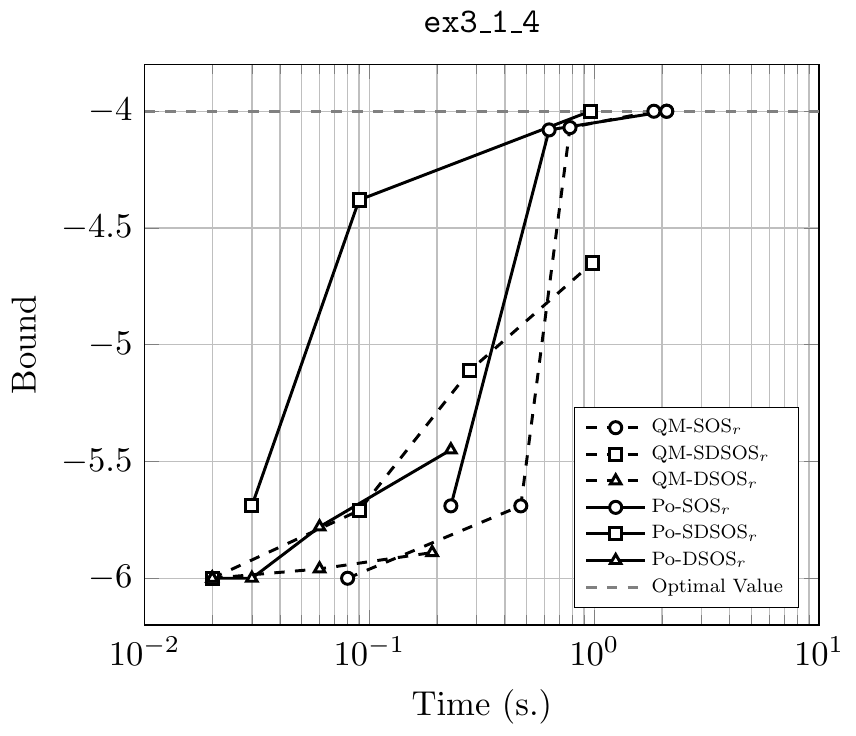}
  %  \caption{Rupture}
    \label{fig7:b}
  \end{subfigure}
  \caption{Bound and Time Comparison of Different Hierarchies for \texttt{ex2\_1\_1} (left) and \texttt{ex3\_1\_4} (right).}
  \label{fig.compare2}
\end{figure}

In Figure~\ref{fig.compare2}, we show the performance of different hierarchies for problem~\verb"ex2_1_1" and problem~\verb"ex3_1_4". Similar to Figure~\ref{fig.compare},
the different line plots show the bound and and solution time associated with increasing orders of each of the hierarchies.
Clearly, for problem~\verb"ex2_1_1", within one second, the~\eqref{opt.sos2} and~\eqref{opt.sdsos2} hierarchies give the optimal value while the bounds obtained by other hierarchies is not tight. Overall, the~\eqref{opt.sdsos2} has the best performance in terms of bound and computational time for problem \verb"ex2_1_1". For problem~\verb"ex3_1_4", within one second, only the~\eqref{opt.sdsos2} reaches the optimal value, again, the~\eqref{opt.sdsos2} has the best performance in terms of bound and computational time for problem \verb"ex3_1_4".

Table~\ref{tab.globlib} shows the bound and time comparison of all hierarchies applied to different test problems. Column~1 shows the name of the problem and column~2 states the number of variables in the problem and its degree, while the degree of each hierarchy $\hat{d}$ is listed in column~3. The results for the Lasserre-type hierarchies are given in columns 4-9 while the remaining columns show the results for the proposed hierarchies in Section~\ref{sec.hierarchy}.
We can see that for problems~\verb"ex_2_1_2",~\verb"ex_2_1_3",~\verb"ex2_1_4", and~\verb"ex2_1_5", the Lassarre-type hierarchies~\eqref{opt.lasserre}, \eqref{opt.socp}, and \eqref{opt.lp} are infeasible when $\hat{d}=2$ and provide the optimal solution when $\hat{d}=4$. In contrast, the proposed hierarchies give the optimal value when $\hat{d}=2$.
\begin{landscape}
\small{
\begin{table}[tb]
    \centering
    \begin{tabular}{cccrrrrrrrrrrrr}
        \toprule
         &  & &\multicolumn{2}{c}{\eqref{opt.lasserre}}&\multicolumn{2}{c}{\eqref{opt.socp}}&\multicolumn{2}{c}{\eqref{opt.lp}}&\multicolumn{2}{c}{\eqref{opt.sos2}}&\multicolumn{2}{c}{\eqref{opt.sdsos2}}&\multicolumn{2}{c}{\eqref{opt.dsos2}}\\
       ID & $(n,d)$ & $\hat{d}$ & Bound & T & Bound & T  & Bound & T & Bound & T & Bound & T & Bound & T\\\midrule
        \verb"ex2_1_1" & (5,2)  &  2 & Infeas. & 0.14 & Infeas.& 0.08 & Infeas. & 0.02 & -18.16 & 0.46 & -18.16 & 0.02 & -18.16 & 0.02 \\
        & & 4 & -17.92 & 0.27 & -18.20 & 0.11 & -18.90 & 0.03 & $^\ast$-17.00 & 0.85 & $^\ast$-17.00 & 0.41 & -17.21 & 0.20   \\
        & & 6 & $^\ast$-17.00 & 1.19 & -17.97 & 0.95 & -18.02 & 0.62 & & &  &  & -17.21 & 0.46   \\\midrule
        \verb"ex2_1_2" & (6,2)  &  2 & Infeas. & 0.11 & Infeas. & 0.73 & Infeas. & 0.00 & $^\ast$-213.00 & 0.64 & $^\ast$-213.00 & 0.03 & $^\ast$-213.00 & 0.02 \\
        & & 4 & $^\ast$-213.00 & 0.46 & $^\ast$-213.00 & 0.25 & $^\ast$-213.00 & 0.06 &   &  &  & &  &  \\ \midrule
        \verb"ex2_1_3" & (13,2)  &  2 & Infeas. & 0.17 & Infeas. & 0.39 & Infeas. & 0.01 & $^\ast$-15.00 & 0.27 & $^\ast$-15.00 & 0.06 & $^\ast$-15.00 & 0.03 \\
        & & 4 & $^\ast$-15.00 & 108.88 & $^\ast$-15.00 & 2.35 & -15.00 & 0.91 &   &  &  & &  &  \\ \midrule
        \verb"ex2_1_4" & (6,2)  &  2 & Infeas. & 0.08 & Infeas. & 0.37 & Infeas. & 0.37 & $^\ast$-11.00 & 0.76 & $^\ast$-11.00 & 0.03 & $^\ast$-11.00 & 0.02 \\
        & & 4 & $^\ast$-11.00 & 0.45 & $^\ast$-11.00 & 0.20 & $^\ast$-11.00 & 0.06 &   &  &  & &  &  \\ \midrule
        \verb"ex2_1_5" & (10,2)  &  2 & Infeas. & 0.13 & Infeas. & 0.75 & Infeas. & 0.01 & $^\ast$-268.01 & 0.26 & $^\ast$-268.01 & 0.08 & $^\ast$-268.01 & 0.03 \\
        & & 4 & $^\ast$-268.01 & 5.05 & $^\ast$-268.01 & 0.86 & -268.06 & 0.84 &   &  &  & &  &  \\ \midrule
        \verb"ex2_1_6" & (10,2)  &  2 & Infeas. & 0.05 & Infeas. & 0.45 & Infeas. & 0.00 & -39.83 & 0.39 & -39.83 & 0.14 & -39.83 & 0.06 \\
        & & 4 & $^\ast$-39.00 & 8.90 & -39.86 & 0.86 & -40.88 & 1.41 & $^\ast$-39.00  & 16.88  & $^\ast$-39.00 & 1.61 & -39.52 & 0.85 \\ \midrule
        \verb"ex2_1_7" & (20,2)  &  2 & Infeas. & 0.28 & Infeas. & 0.36 & Infeas. & 0.02 & -4331.16 & 1.10 & -4332.53 & 0.61 & -4332.53 & 0.16 \\
         & &4 & $\diamond$  &  $\diamond$  &-4686.70  & 23.08 &$-\infty$  &8.79  &$\circ$  &$\circ$  &$\circ$  &$\circ$  &$\circ$&$\circ$  \\ \midrule
        \verb"ex2_1_10" & (20,2)  &  2 & Infeas. & 0.22 & Infeas. & 0.38 & Infeas. & 0.02 & $^\ast$49318.02 & 0.69 & $^\ast$49318.02 & 0.25 & $^\ast$49318.02 & 0.06 \\
        & & 4 & $\diamond$ & $\diamond$& $^\ast$49318.02 & 22.92 & $-\infty$ & 7.70 &  &  &  &  & &  \\\midrule
        \verb"ex3_1_3" & (6,2) &  2 & Infeas. & 0.11 & Infeas.& 0.74 & Infeas. & 0.00 & $^\ast$-310.00 & 0.67 & $^\ast$-310.00 & 0.02 & $^\ast$-310.00 & 0.02 \\
        & & 4 & $^\ast$-310.00 & 0.71 & $^\ast$-310.00 & 0.38 & $-\infty$ & 0.08 &  &  &  & &  &  \\ \midrule
        \verb"ex3_1_4" & (3,2)  &  2 & -6.00 & 0.08 & -6.00 & 0.02 & -6.00 & 0.02 & -5.69 & 0.23 & -5.69 & 0.03 & -6.00 & 0.02 \\
        & & 4 & -5.69 & 0.47 & -5.71 & 0.09 & -6.00 & 0.02 & -4.08 & 0.63 & -4.38 & 0.09 & -6.00 & 0.03  \\
        & & 6 & -4.07 & 0.78 & -5.11 & 0.28 & -5.96 & 0.06 & $^\ast$-4.00  & 2.10 & $^\ast$-4.00 & 0.96 & -5.78 & 0.06   \\
        & & 8 & $^\ast$-4.00 & 1.84 & -4.65 & 0.98 & -5.89 & 0.19 &   &  &  &  & -5.45 & 0.23   \\\midrule
        \verb"ex4_1_9" & (2,4) &  4 & -7.00 & 0.15 & -7.00 & 0.04 & -7.00 & 0.02 & $^\ast$-5.51 & 0.20 & $^\ast$-5.51 & 0.04 & -7.00 & 0.02 \\
        & & 6 & $^\ast$-5.51 & 0.23 & $^\ast$-5.51 & 0.07 & -7.00 & 0.02 &   &  &  & & $\circ$  & $\circ$ \\ \midrule
        \verb"mathopt2" & (2,4)  &  4 & Infeas. & 0.45 & Infeas. & 0.71 & Infeas. & 0.13 & $^\ast$0.00 & 0.21 & $^\ast$0.00 & 0.03 & $^\ast$0.00 & 0.02 \\
        & & 6 & $^\ast$0.00 &  0.17 & $^\ast$0.00 & 0.05 & $-\infty$ & 0.02 &   &  &  & &  &  \\
        \bottomrule
    \end{tabular}
    \begin{tablenotes}
        \footnotesize
        \item[1] $\ast$: Optimal value is obtained.
        \item[2] $\diamond$: Solver runs out of memory.
        \item[3] $\circ$: {\tt Matlab} runs out of memory while formulating LMI.
    \end{tablenotes}
\caption{Bound and Time Comparison of Different Hierarchies for Examples in Global Optimization Library.}\label{tab.globlib}
\end{table}
}
\end{landscape}
The time to obtain the optimal value is greatly reduced by using~\eqref{opt.sos2} over~\eqref{opt.lasserre} for problem \verb"ex_2_1_3". For problems~\verb"ex_2_1_7" and~\verb"ex2_1_10" with relatively large numbers of variables (20 variables), the Lasserre-type hierarchies are infeasible when $\hat{d}=2$, which means that the Lasserre-type hierarchies fail to give a bound, however, by using the hierarchies proposed here, the optimal value for problem \verb"ex2_1_10" and a global lower bound for problem~\verb"ex_2_1_7" are obtained when $\hat{d}=2$. The hierarchy~\eqref{opt.lasserre} runs out of memory when $\hat{d}=4$ for problems~\verb"ex_2_1_7" and \verb"ex2_1_10". For~\verb"ex_2_1_7", the hierarchy~\eqref{opt.socp} gives a bound when~$\hat{d}=4$; however, it is weaker than the ones obtained from the hierarchies~\eqref{opt.sos2},~\eqref{opt.sdsos2} and~\eqref{opt.dsos2} when $\hat{d}=2$. {\tt Matlab} runs out of memory when formulating the LMI for the hierarchies~\eqref{opt.sos2},~\eqref{opt.sdsos2} and~\eqref{opt.dsos2} when~$\hat{d}=4$ for~\verb"ex_2_1_7", due to a large number of constraints in~\verb"ex_2_1_7". For other cases, our proposed hierarchies mostly converge to global optimum with a smaller~$\hat{d}$ than that of Lasserre-type hierarchies.

From Table~\ref{tab.globlib}, one can notice that in some instances of the problems, there is no improvement in the bound obtained by using the sequentially tighter~\eqref{opt.sos2},~\eqref{opt.sdsos2} and~\eqref{opt.dsos2} hierarchies. For some problems (like {\tt ex2\_1\_2}), this is due to the~\eqref{opt.dsos2} hierarchy providing the problem's optimal solution. For other problems (like {\tt ex2\_1\_1}), this is a result of the structure of the problem, which results in hierarchies~\eqref{opt.sos2},~\eqref{opt.sdsos2} not helping to improve the bounds. For example, for problem {\tt ex2\_1\_1}, the objective function is given by~$f(x_1,\dots,x_5)=42x_1+44x_2+45x_3+47x_4+47.5x_5-50(x_1^2+x_2^2+x_2^2+x_4^2+x_5^2)$. Since there are no cross-variable terms, the bound obtained from the hierarchies \eqref{opt.sos2}, \eqref{opt.sdsos2} does not take advantage of providing a tighter formulation of the POP for cross-variable monomials. On the other hand, it is clear that in problems like {\tt ex3\_1\_4}, the tighter  \eqref{opt.sos2}, \eqref{opt.sdsos2} hierarchies provide better bounds than the \eqref{opt.dsos2} hierarchy.

\subsection{Numerical Results on Problems with More Variables}
Consider the following non-convex problem,
\bequation\label{opt.numexp}
    \baligned
              \min_{} & &   & \sum_{|\alpha|\leq 2} c_\alpha x^\alpha \\
                 \st &  &  &\sum_{i=1}^nx_i^2\leq 1,\\
                      & & & \sum_{i=1}^nx_i^2\geq 0.6^2,\\
                    &    & & x_i\geq 0, i =1,\ldots,n,
    \ealigned
\eequation
where $c_\alpha$ are randomly generated in $[-1,1]$. We construct this nonconvex problem inspired by~\cite{yang2016quadratic} to test the performance of the proposed hierarchies. The problem is to find the minimal value of a polynomial on the Euclidean unit ball intersected with the positive orthant while excluding the Euclidean ball with radius 0.6. We use instances with relatively large number of variables $n\in[20,30,50,100,150]$ and compare the Lasserre-type hierarchies with the hierarchies proposed in Section~\ref{sec.hierarchy}. Note that it will be computationally expensive to run higher levels of the hierarchies for large-scale problems. The purpose of this comparison mainly focuses on the bound obtained for quadratic programs when $\hat{d}=2$.

\begin{figure}[ht]
    \centering
    \includegraphics[width=0.5\linewidth]{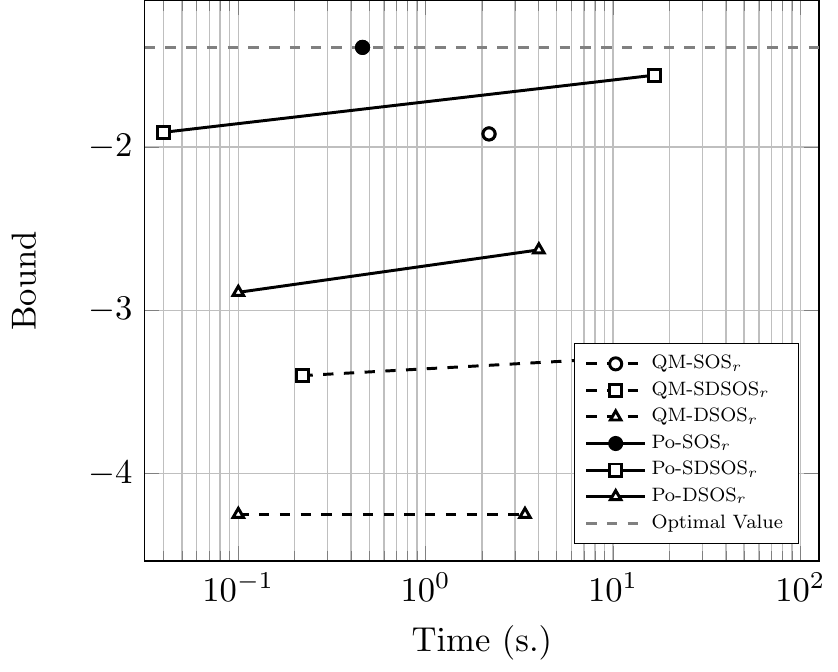}
  %  \caption{Rupture}
  \caption{Bound and Time Comparison of Different Hierarchies for Problem (7) with $n=20$.}
  \label{fig.compare3}
\end{figure}

Similar to the figures in previous sections, in Figure~\ref{fig.compare3}, we plot the performance of all hierarchies for problem~\eqref{opt.numexp} with $n=20$. Clearly, only the~\eqref{opt.sos2} obtains the optimal value with $n=20$. Unlike the instances in Section~\ref{sec.numsub1} and Section~\ref{sec.numsub2}, the increasing order of the hierarchies doesn't improve the bound significantly, thus the lines in~Figure~\ref{fig.compare3} are flat.

Table~\ref{tab.numexp} lists the bound and time of all hierarchies for problem~\eqref{opt.numexp} with different $n$. The results for $n\geq 30$ and $\hat{d}=4$ are not listed since {\tt Matlab} runs out of memory when formulating the LMI for these cases. Column 2 is the upper bound we obtain from a global optimization solver, BARON \cite{baron}; columns 4-15 list the results by running Lasserre-type hierarchies and the proposed hierarchies. The optimal value (indicated by $\ast$) is obtained when the upper bound obtained by BARON is equal to any of the lower bound from the six hierarchies. It is clear that our proposed hierarchies~\eqref{opt.sos2},~\eqref{opt.sdsos2}, and~\eqref{opt.dsos2} yield tighter bounds than corresponding the Lasserre-type hierarchies. For cases with $n=20,30,50,100$, the hierarchy~\eqref{opt.sos2} converges to the global optimum when $\hat{d}=2$. For the case with $n=150$, SOS-based hierarchies~\eqref{opt.lasserre} and~\eqref{opt.sos2} fail to give a bound due to the computationally difficult SDP constraints. SOCP-based hierarchies can be used to obtain global bounds, in which case our proposed hierarchy~\eqref{opt.sdsos2} improves the bounds obtained from~\eqref{opt.socp} by approximately 100\%. LP-based hierarchies provide the worst bounds among the same type of hierarchies, however, the bound obtained by the LP-based hierarchy~\eqref{opt.dsos2} is even tighter than the bound obtained by SOCP-based hierarchy~\eqref{opt.socp}.

\begin{table}[ht]
    \centering
   \begin{adjustbox}{width=1\textwidth}
    \begin{tabular}{cccrrrrrrrrrrrr}
        \toprule
        & BARON &  &\multicolumn{2}{c}{\eqref{opt.lasserre}}&\multicolumn{2}{c}{\eqref{opt.socp}}&\multicolumn{2}{c}{\eqref{opt.lp}}&\multicolumn{2}{c}{\eqref{opt.sos2}}&\multicolumn{2}{c}{\eqref{opt.sdsos2}}&\multicolumn{2}{c}{\eqref{opt.dsos2}}\\
        $(n,d)$  & UB & $\hat{d}$  & LB & T & LB & T  & LB & T & LB & T & LB & T & LB & T\\\midrule
        (20,2)  &  $^\ast$-1.39 & 2  &  -1.92 & 2.18 & -3.40 & 0.22 & -4.25 & 0.01 & $^\ast$-1.39 & 0.46 & -1.91 & 0.04 & -2.89 & 0.01 \\
         & & 4 & $\diamond$ & $\diamond$ & -3.29 & 11.99 & -4.25 & 3.39 &  &  & -1.56 & 16.67 & -2.63 & 4.02  \\\midrule
         (30,2) & $^\ast$-1.96 &  2 & -2.21 & 2.61 & -4.82 & 0.16 & - 6.40 & 0.01 & $^\ast$-1.96 & 2.08 & -2.91 & 0.06 & -4.93 & 0.01 \\ \midrule
          (50,2)  & $^\ast$-2.06 & 2 & -2.90 & 43.99 & -8.30 & 1.82 & -11.38 & 0.38 & $^\ast$-2.06 & 50.26 & -4.32 & 0.15 & -6.38 & 0.01 \\ \midrule
           (100,2)  & $^\ast$-3.02 &  2 & -3.89 & 1606.00 & -16.56 & 1.12 & -20.47 & 0.01 & $^\ast$-3.02 & 1824.30 & -8.49 & 0.45 & -11.92 & 0.04 \\ \midrule
         (150,2)  & -3.58 &  2 & $\diamond$ & $\diamond$ & -25.02 & 0.90 & -29.12 & 0.04 & $\diamond$ & $\diamond$ & -12.66 & 1.11 & -16.42 & 0.09 \\
        \bottomrule
    \end{tabular}
  \end{adjustbox}
    \begin{tablenotes}
        \footnotesize
        \item[1] LB,UB: lower bound, upper bound.
         \item[2] $\ast$: Optimal value is obtained.
        \item[3] $\diamond$: Solver runs out of memory.
    \end{tablenotes}
\caption{Bound and Time Comparison of Different Hierarchies for Problem~\eqref{opt.numexp}.}\label{tab.numexp}
\end{table}

%*********
\section{Concluding Remarks}\label{sec.conclusion}
In this paper, we propose alternative LP, SOCP and SDP approximation hierarchies to obtain global bounds for general POPs, by using SOS, SDSOS and DSOS polynomials to strengthen existing LP hierarchy for POPs. Comparing with the classic Lasserre's hierarchy, the LP and SOCP approximation hierarchies are shown to be computationally more efficient to find the global optimum of POPs for which Lasserre's hierarchy is known to perform poorly. In particular, this shows that the relaxation approach introduced by~\citeauthor{ref.amirali}~\cite{ref.amirali} produces better results as a way to strengthen LP-based hierarchies for POPs. Furthermore, these hierarchies are shown to converge as the level of the hierarchy increases to the global optimum of the corresponding POP. Unlike other hierarchies proposed in the literature, this property is obtained whenever the feasible set of the POP is compact but the {\em quadratic module} of the polynomials defining the problem's feasible set is not necessarily{\em Archimedean}.

The fact that the hierarchies considered here are based on using LP and SOCP allows for the future use of column
generation approaches in order to be able to address the solution of larger-scale POPs.

%**************
% Bibliography
%**************
\bibliographystyle{apalike}
\bibliography{references}

\end{document}